\documentclass[12pt]{article}
\usepackage{amssymb,amsmath, amsthm}
\usepackage{graphicx}

\newtheorem{thm}{Theorem}

\newtheorem*{propo}{Proposition}

\newtheorem*{coro}{Corollary}
\newtheorem{lemma}{Lemma}
\newtheorem*{lem}{Lemma}
\newenvironment{defin}{\medskip\noindent{\sc
Definition}. }{\goodbreak\medskip}
\newenvironment{nota}{\medskip\noindent{\sc
Notations}. }{\goodbreak\medskip}
\newenvironment{remk}{\noindent{\sc
Remark}. }{\goodbreak\vskip10pt}
\newtheorem{prop}{Proposition}

\newtheorem{ques}{Question}

\def\demo{\medskip\goodbreak\noindent
     \hbox{\sc Proof \kern .3em}\ignorespaces}%
  \def \qedbox{$\square$}%
  \def \qed{\hglue1mm\hfill{\ifmmode\qedbox
     \else\unskip\ \hglue0mm\hfill\qedbox\medskip
      \goodbreak\fi}}%
\def\enddemo{\qed\goodbreak\vskip10pt}%

\def\qed{\hglue1mm\hfill\raise -2pt\hbox{\vrule\vbox to 10pt{\hrule width
4pt
                  \vfill\hrule}\vrule}}

\newcommand{\R}{\mathbb {R}}

\newcommand{\N}{\mathbb {N}}

\newcommand{\Hc}{\mathcal {H}}

\newcommand{\Ec}{\mathcal {E}}

\newcommand{\Lc}{\mathcal {L}}

\newcommand{\Ac}{\mathcal {A}}

\begin{document}
\title{Rigidity in topology  $C^0$ of the Poisson bracket for Tonelli Hamiltonians}
\author{M.-C. ARNAUD 
\thanks{ANR-12-BLAN-WKBHJ}
\thanks{Avignon Universit\'e , Laboratoire de Math\'ematiques d'Avignon (EA 2151),  F-84 018 Avignon,
France. e-mail: Marie-Claude.Arnaud@univ-avignon.fr} 
\thanks{membre de l'Institut universitaire de France}}
\maketitle
\abstract{\noindent \sl We prove the following rigidity result for the Tonelli Hamiltonians.\\
Let $T^*M$ be the cotangent bundle of a closed manifold $M$ endowed with its usual symplectic form. Let $(F_n)$ be a sequence of Tonelli Hamiltonians that $C^0$ converges on the compact subsets to a Tonelli Hamiltonian $F$. Let $(G_n)$ be a sequence of Hamiltonians that that $C^0$ converges on the compact subsets to a Hamiltonian $G$. We assume that the sequence of the Poisson brackets $(\{ F_n, G_n\})$ $C^0$-converges on the compact subsets to a $C^1$ function $\Hc$. Then $\Hc=\{ F, G\}$.}

\newpage
   \section{Introduction}  For two $C^2$ functions $H,K: N\rightarrow \R$ defined on a symplectic manifold $(N, \omega)$ (such functions are usually called ``Hamiltonians''), their Poisson bracket $\{ H, K\}$ is  
 $$\{ H, K\} =\omega (X_H, X_K)$$ where $X_H$ designates the symplectic gradient of $H$ defined by $dH=\omega(X_H, .)$.  Hence the Poisson bracket is given by $\{ H, K\} =dH(X_K)$ and describes the evolution of $H$ along any orbit for the Hamiltonian flow of $K$. The Poisson bracket is  clearly continuous for the topologies $(H, K)\in C^1(N, \R)\times C^1( N, \R)\mapsto \{ H, K\}\in C^0(N, \R)$.
 
 In his PhD thesis \cite{Hum1}, V.~Humil\`ere asks  the following question concerning a possible $C^0$ rigidity for the Poisson bracket: 
 
\begin{ques}\label{Qrigidity}{ (V.~Humili\`ere) Let $(F_k)$, $(G_k)$ be two sequences of $C^\infty$ functions such that $(F_k)$ $C^0$-converges  to $F\in C^\infty(N, \R)$, $(G_k)$ $C^0$-converges  to $G\in C^\infty(N, \R)$ and $(\{F_k, G_k\})$ $C^0$-converges  to $\Hc\in C^\infty(N, \R)$. Do we have $\{ F, G\}=\Hc$?}
 \end{ques}

 The answer to this question is in general negative and some counter-examples due to L.~Polterovitch are provided in \cite{Hum2}.\\
  However, F.~Cardin and C.~Viterbo give a positive answer in \cite{CarVit} when $\Hc=0$ and all the considered Hamiltonians have their support contained in fixed  a compact set.\\
  By introducing a notion of pseudo-representation, V.~Humili\`ere obtains  in \cite{Hum2} results that contain the previous one. Fix a closed manifold $M$ and a finite-dimensional Lie algebra $({\bf g}, [.,.])$, a pseudo-representation is a sequence  of morphisms $\rho_n:{\bf g}\rightarrow C^\infty(M, \R)$ such that for all $f, g\in{\bf g}$, the sequence $(\{\rho_n(f), \rho_n(g)\}-\rho_n([f, g]))$ converges uniformly to zero. Humil\`ere proves that the limit of a convergent pseudo-representation is a representation. In other word, assuming that there is a morphism $\rho: {\bf g}\rightarrow C^\infty(M, \R)$ such that any sequence $(\rho_n(f))$ uniformly converges to $\rho(f)$, he proves that
$$\forall f, g\in {\bf g}, \{ \rho(f), \rho(g)\}=\rho([f, g]);$$
i.e. that $\displaystyle{  \{\lim_{n\rightarrow \infty} \rho_n(f), \lim_{n\rightarrow \infty} \rho_n(g)\}=\lim_{n\rightarrow \infty}\{\rho_n(f), \rho_n(g)\}}$.

      In \cite{Buh1},  L.~Buhovski proves that if $(f_n)$ and $(g_n)$ are $C^\infty$ functions that uniformly converge on every compact subset of $N$ to the smooth functions $f$ and $g$, and if the sequences $(Df_n)$, $(Dg_n)$ are uniformly bounded on any compact subset of $N$, then we have 
      $$\lim_{n\rightarrow \infty}\{ f_n, g_n\}=\{ f, g\}.$$      \\
 We will   describe here a new   case where question \ref{Qrigidity} always has a positive question: the case where one of the considered sequence is a sequence of Tonelli Hamiltonians that $C^0$-converges to a Tonelli Hamiltonian. We will too give a simpler proof of a part of Buhovski's result.
   
 \nota Let $M$ be a closed manifold and let $T^*M$ be its cotangent bundle endowed with its usual symplectic form $\omega$ that is the derivative of the Liouville 1-form.  We use the notation $(q, p)\in T^*M$ and then $\omega=dq\wedge dp$. \\
 The first projection is denoted by $\pi:(q, p)\in T^*M\mapsto q\in M$.
 
 \begin{defin} A $C^2$ function $H: T^*M\rightarrow \R$ is   a {\em Tonelli Hamiltonian} if it is: 
\begin{itemize} 
\item superlinear in the fiber, i.e. $\forall A\in \R, \exists B\in \R, \forall (q,p)\in T^*M, \| p\|\geq B\Rightarrow  H(q,p) \geq A\| p\|$;
\item $C^2$-convex in the fiber i.e. for every $(q,p)\in T^*M$, the Hessian $\frac{\partial^2H}{\partial p^2}$ of $H$ in the fiber direction is positive definite as a quadratic form.
\end{itemize}
We denote the Hamiltonian flow of $H$ by $(\varphi^H_t)$ and the Hamiltonian vector-field by $X_H$.
 \end{defin}
 Note that the class of Tonelli Hamiltonian contains all the Riemannian metrics and all the mechanical systems (Riemannian metric + potential). The class of Tonelli Hamiltonians has been recently intensely studied in the setting of weak KAM theory.
 
 \begin{thm}\label{Tonelli}
 Let $(F_k)$  be a sequence of Tonelli Hamitonians defined on $T^*M$  such that $(F_k)$ $C^0$-converges    on the compact subsets of $T^*M$ to a Tonelli Hamiltonian $F\in C^2(T^*M, \R)$, let $(G_k)$ be a sequence of $C^2$ Hamiltonians of $T^*M$ that  $C^0$-converges     on the compact subsets of $T^*M$  to $G\in C^2(T^*M, \R)$ and such that $(\{F_k, G_k\})$  $C^0$-converges      on the compact subsets of $T^*M$ to $\Hc\in C^0(T^*M, \R)$. Then  $\{ F, G\}=\Hc$.
 \end{thm}
\remk We prove in lemma \ref{Roc} that the $C^0$ convergence of $(F_k)$ to $F$ and the convexity hypothesis imply the $C^0$-convergence of $\frac{\partial F_k}{\partial p}$ to $\frac{\partial F}{\partial p}$. But we don't have the convergence of the $\frac{\partial F_k}{\partial q}$; then this theorem is not a consequence of theorem \ref{C1}.\medskip

 \begin{thm} \label{C1}
 Let $(N, \omega)$ be a (non necessarily  compact) symplectic  manifold. Let $(F_k)$, $(G_k)$ be two sequences of $C^2$ Hamitonians defined on $N$  such that $(G_k)$ $C^0$-converges     on the compact subsets of $N$ to $G\in C^2(N, \R)$, $(F_k)$ $C^1$-converges     on the compact subsets of $N$  to $F\in C^2(N, \R)$ and $(\{F_k, G_k\})$ $C^0$-converges  on the compact subsets of $N$ to $\Hc\in C^0(N, \R)$. Then  $\{ F, G\}=\Hc$.
 \end{thm} 
 Theorem 25.7 of \cite{Roc1}� asserts that if $(F_k)$ is a sequence of $C^1$ and convex functions defined on $\R^n$ that $C^0$-converge to a $C^1$ and convex function $F$ on any compact subset, then the convergence is $C^1$ on any compact subset of $\R^n$. We deduce:
 \begin{coro}
 Let $(F_k)$, $(G_k)$ be two sequences of $C^2$ Hamitonians defined on $\R^{2n}$  such that $(G_k)$ $C^0$-converges     on the compact subsets of $\R^{2n}$ to $G\in C^2(\R^{2n}, \R)$, $(F_k)$ $C^0$-converges     on the compact subsets of $\R^{2n}$  to $F\in C^2(\R^{2n}, \R)$, the $F_k$ are convex  and $(\{F_k, G_k\})$ $C^0$-converges  on the compact subsets of $\R^{2n}$ to $\Hc\in C^0(\R^{2n}, \R)$. Then  $\{ F, G\}=\Hc$.
 \end{coro} 

Theorem \ref{Tonelli} concerns Tonelli Hamiltonians and is new. Theorem \ref{C1} is a consequence of the result due to L.~Buhovsky (see \cite{Buh1}) that we explained before, but we give here a simple proof, that uses the simple principle.
 \subsection*{Ideas of proofs}
 Contrarily to what happens in the works that we mentioned above, we won't use here any  result of symplectic topology, as symplectic capacities, displacement energy\dots \\
  We use
  \begin{enumerate}
  \item   the following simple principle:\\
  
\noindent  {\em {\bf Simple principle:} under the hypotheses of question \ref{Qrigidity}, if for any $x\in N$, there exists a $T>0$ and a sequence of points $(x_k)$ such that the arcs of orbits $(\varphi^{F_k}_t(x_k))_{t\in [0, T]}$ converge in some reasonable sense to $(\varphi^{F}_t(x))_{t\in [0, T]}$, then the answer to question \ref{Qrigidity} is positive. }\\
If we have a good enough notion of convergence, the proof of this principle is not  complicated and propositions \ref{clef} and \ref{2clef} are some versions of this principle.
\item for theorem \ref{Tonelli}, we use some subtle variational arguments to prove that the hypothesis of the simple principle is satisfied; for theorem \ref{C1}, we just use straightforward arguments coming from the theory of ordinary differential equations.
\end{enumerate}

\subsection*{Plan}
\begin{enumerate}
\item[$\bullet$] In section \ref{Sec1}, we give the precise statements of the propositions that we need to prove theorems \ref{Tonelli} and \ref{C1}.\\
\item[$\bullet$] Section \ref{Sec3} is devoted to the proof of proposition \ref{Pconver} that asserts  that every $C^0$-convergent sequence of Tonelli Hamiltonians satifies the hypotheses of the simple principle.\\
\item[$\bullet$] In section \ref{Sec4}, we prove the simple principle, i.e.  propositions \ref{clef} and \ref{2clef}, in the considered cases. 
\end{enumerate}
{\bf Acknowledgments.} I am grateful   to the referees for many improvements of
the article.

 \section{Structure of the proofs of theorems \ref{Tonelli} and \ref{C1}}\label{Sec1}
 
\subsection{ How two propositions imply theorem \ref{Tonelli}}
We will prove the two following propositions in sections \ref{Sec3} and \ref{Sec4}. The first one uses variational arguments to state that if a sequence of Tonelli Hamiltonians $(H_k)$ $C^0$-converges to a Tonelli Hamiltonian $H$, then the ``small'' arcs of orbits for $H$ are approximated in some sense by some small parts of orbits for the $H_k$.

 \begin{prop}\label{Pconver}
 Let $(H_k)_{k\in\N}$, $H$ be some $C^2$ Tonelli Hamiltonians defined on $T^*M$ such that $(H_k)$ $C^0$-converges   on the compact subsets of $T^*M$ to $H$. 
 
 Then, for every $(q,p)\in T^*M$, there exists $T>0$ and $(q_k, p_k)_{k\in \N}\in T^*M$  such that a subsequence of the  sequence of arcs $(Q_k, P_k)=(\varphi_t^{H_k}(q_k, p_k))_{t\in [0, T]}$ converges to  $(Q, P)=(\varphi_t^H(q, p))_{t\in[0, T]}$ in the following sense (here we are in a chart because all is local and we denote the subsequence in the same way we denoted the initial sequence):
 \begin{enumerate}
 \item[$\bullet$] $(Q_k)$ $C^0$-converges to $Q$;
 \item[$\bullet$] $(P_k)$ converges to $P$ in the $L^2$ sense;
 \item[$\bullet$] the $(Q_k, P_k)$ are uniformly $C^0$ bounded.
 \end{enumerate} 
 \end{prop} \label{piece}
 
 The second proposition is a  version of what we called previously the {\em simple principle}. The proof is elementary.
 
  \begin{prop}\label{clef}  Let $(F_k)$ be a sequence of $C^2$ Hamiltonians defined on $T^*M$ and let $F$ be a $C^2$ Hamiltonian such that:\\
   for every $(q,p)\in T^*M$, there exists $T>0$ and $(q_k, p_k)_{k\in \N}\in T^*M$  such that the sequence of arcs $(Q_k, P_k)=(\varphi_t^{F_k}(q_k, p_k))_{t\in [0, T]}$ converges to  $(Q, P)=(\varphi_t^F(q, p))_{t\in[0, T]}$ in the following sense:
 \begin{enumerate}
 \item[$\bullet$] $(Q_k)$ $C^0$-converges to $Q$;
 \item[$\bullet$] $(P_k)$ converges to $P$ in the $L^2$ sense;
 \item[$\bullet$] the $(Q_k, P_k)$ are uniformly $C^0$ bounded.
 \end{enumerate} 
Let $(G_k)$ be a sequence of $C^2$ Hamiltonians of $T^*M$ that $C^0$-converges on the compact subsets of $T^*M$ to a $C^2$ Hamiltonian $G$ and that is such that the sequence $(\{ F_k, G_k\})$ $C^0$-converges to some $\Hc\in C^0(T^*M, \R)$.\\
Then $\{ F, G\} =\Hc$.

 \end{prop} 
 Theorem \ref{Tonelli} is clearly a consequence of propositions \ref{Pconver} and \ref{clef}.
 
 \subsection{How two propositions imply theorem \ref{C1}}
 To prove theorem \ref{C1}, we will prove that if $(F_k)$ satisfy the hypothesis of theorem \ref{C1}, then it satisfies the hypothesis of proposition \ref{2clef} below .

  \begin{prop}\label{2clef} Let $(N, \omega)$ be a (non necessarily compact) symplectic manifold. Let $(F_k)$ be a sequence of $C^2$ Hamiltonians defined on $N$ and let $F$ be a $C^2$ Hamiltonian such that:\\
   for all $x\in N$, there exists  $T>0$ and a sequence $(x_k)\in N$ such that the sequence of arcs of orbit $(\varphi_t^{F_k}(x_k))_{t\in [0, T]}$ $C^0$-converges to $(\varphi_t^{F}(x))_{t\in [0, T]}$.\\
Let $(G_k)$ be a sequence of $C^2$ Hamiltonians of $N$ that $C^0$-converges on the compact subsets of $N$ to a $C^2$ Hamiltonian $G$ and that is such that the sequence $(\{ F_k, G_k\})$ $C^0$-converges on the compact subsets of $N$ to some $\Hc\in C^0(N, \R)$.\\
Then $\{ F, G\} =\Hc$.

 \end{prop} 
 Let us be more precise.
 
 \begin{defin}
 Let $X$ be a locally Lipschitz vectorfield on $N$. A $C^1$ arc $\gamma: [a, b]\rightarrow N$ is an $\varepsilon$-solution for $X$ if $\forall t\in [a,  b], \| \gamma'(t)-X(\gamma(t))\|<\varepsilon$.
  \end{defin}
  If for example the vector fields $X_1$ and $X_2$ are $\varepsilon$-close, every solution for $X_2$ is an $\varepsilon$-solution for $X_1$.
  
  The following proposition is classical (see for example section 6.2.2. of  \cite{HubWes}):
  
 \begin{prop}\label{C1orb} Assume that the vector field $X$ is  Lipschitz on $\R^n$ with Lipschitz constant $K$.  Let $\gamma_i: [a, b]\rightarrow \R^n$ be a $\varepsilon_i$-solution for $X$  for $i=1, 2$. We assume that for some $\tau\in [a, b]$, we have $d(\gamma_1(\tau), \gamma_2(\tau))\leq\delta$. Then we have:
 $$\forall t\in [a, b], \quad d(\gamma_1(t), \gamma_2(t))\leq  \delta e^{K|t-\tau|}+\frac{\varepsilon_1+\varepsilon_2}{2}(e^{K|t-\tau|}-1).$$

 \end{prop}

 Let now $(F_k)$, $F$ be $C^2$ Hamiltonians on $N$ such that $(F_k)$ $C^1$-converges  to $F$. Using proposition \ref{C1orb}, we will prove:\\
\begin{lem} for all $x\in N$, there exists  $T>0$ and a sequence $(x_k)\in N$ such that the sequence of arcs of orbit $(\varphi_t^{F_k}(x_k))_{t\in [0, T]}$ $C^0$-converges to $(\varphi_t^{F}(x))_{t\in [0, T]}$. 
\end{lem}
 Then the conclusion of theorem \ref{C1} follows from this and proposition \ref{2clef}.

 \noindent{\bf Proof of the lemma.}
Let us fix $x_0\in N$. Then we choose a chart $U=B(x_0, r)$ at $x_0$, and $T>0$ such that $\displaystyle{ \sup_{x\in B(x_0, r)}\| X_F(x)\|.T <r}$. As $(F_k)$ $C^1$-converges to $F$, we have  for $k$ large enough: $\displaystyle{ \sup_{x\in B(x_0, r)}\| X_{F_k}(x)\|.T <r}$. This implies that $$\forall t\in [0, T], \varphi_t^{F_k}(x_0)\in B(x_0, r).$$
Moreover, $(\varphi_{t}^{F_k}(x_0))_{t\in [0, T]}$ is a $\varepsilon_k$-solution for $X_{F_k}$ where $\displaystyle{\varepsilon_k=\sup_{x\in B(x_0, r)}\| X_{F_k}(x)-X_F(x)\|}$. We deduce from proposition \ref{C1orb} that for $k$ large enough
$$\forall t\in [0, T], d(\varphi_t^{F_k}(x_0), \varphi_t^{F}(x_0))\leq \frac{\varepsilon_k}{2}(e^{Kt}-1) $$
where $K$ is a Lipschitz constant for $X_{F|B(x_0, r)}$. As $\displaystyle{\lim_{k\rightarrow \infty} \varepsilon_k=0}$, we obtain the wanted conclusion.
\enddemo

 \section{Proof of proposition \ref{Pconver}}\label{Sec3}
  
 We assume that $(H_k)$ $C^0$-converges   on the compact subsets of $T^*M$ to $H$.  We will use the following lemma, that is in the spirit of theorem 25.7 of \cite{Roc1}.
 \begin{lemma}\label{Roc} Let $(F_k)$ be a sequence of $C^2$ functions defined on $T^*M$ that are $C^2$ convex  in the fiber direction and that $C^0$-converges on the compact subsets of $T^*M$ to a $C^2$ function $F:T^*M\rightarrow \R$ that is $C^2$ convex in the fiber direction. Then $(\frac{ \partial F_k}{\partial p} )$ $C^0$-converges   to $\frac{\partial F}{\partial p}$ on the compact subsets of $T^*M$.
 
 \end{lemma}
 \noindent{\bf Proof of lemma \ref{Roc}.} If not, by possibly extracting  a subsequence,  we find a compact subset $K\subset T^*M$, $\varepsilon >0$ and a sequence $(q_k, p_k)\in K$ such that:
 $$\forall k; \| \frac{\partial F_k}{\partial p}(q_k, p_k)-\frac{\partial F}{\partial p}(q_k, p_k)\|\geq 2\varepsilon.$$
 Extracting a subsequence, we can assume that $\displaystyle{ \lim_{k\rightarrow \infty}(q_k, p_k)=(q_\infty, p_\infty)}$ and that we are in a chart $K=K_0\times \bar B(p_\infty, R)\subset M\times \R^n$.
 Then we choose $\tilde K= \tilde K_0\times \bar B(p_\infty, 2R)$ that contains $K$ in its interior. We choose $\alpha>0$ such that $\bar B_\alpha=\bar B((q_\infty, p_\infty), 2\alpha)\subset \tilde K$ and:
 $$\forall x\in \bar B_\alpha, \left\| \frac{\partial F}{\partial p}(x)-\frac{\partial F}{\partial p}(q_\infty, p_\infty)\right\| <\eta=\frac{\varepsilon}{3}.$$
 As $ \| \frac{\partial F_k}{\partial p}(q_k, p_k)-\frac{\partial F}{\partial p}(q_k, p_k)\|\geq 2\varepsilon$, we can find $u_k\in\R^n$ such that $\| u_k\|=1$ and 
 $$\left(\frac{\partial F_k}{\partial p}(q_k, p_k)-\frac{\partial F}{\partial p}(q_k, p_k)\right).u_k\geq \varepsilon;$$
 i.e.
 $$\frac{\partial F_k}{\partial p}(q_k, p_k).u_k \geq \frac{\partial F}{\partial p}(q_k, p_k).u_k +\varepsilon.$$
 Using the fact that $F_k$ is $C^2$ convex in the fiber direction, we deduce:
 $$\begin{matrix}F_k(q_k, p_k+\alpha.u_k)&\geq & F_k(q_k, p_k)+\alpha\frac{\partial F_k}{\partial p}(q_k, p_k).u_k\\
 &\geq & \alpha.\varepsilon + F_k(q_k, p_k)+\alpha \frac{\partial F}{\partial p}(q_k, p_k).u_k\\
 &\geq & \alpha\varepsilon +  F_k(q_k, p_k)+ \alpha \frac{\partial F}{\partial p}(q_\infty, p_\infty).u_k-\alpha\eta
 \end{matrix}$$
 Extracting a subsequence, we can assume that $(u_k)$ converges to $u_\infty$. We then take the limit and obtain:
 $$F(q_\infty, p_\infty+\alpha.u_\infty)- F(q_\infty, p_\infty)-\alpha \frac{\partial F}{\partial p}(q_\infty, p_\infty)u_\infty\geq \alpha(\varepsilon-\eta)=\frac{2\alpha\varepsilon}{3}.$$
 But we have:
 $$\begin{matrix} F(q_\infty, p_\infty+\alpha.u_\infty)&-F(q_\infty, p_\infty)-\alpha \frac{\partial F}{\partial p}(q_\infty, p_\infty)u_\infty\hfill\\
 &=\alpha.\int_0^1\left( \frac{\partial F}{\partial p}(q_\infty, p_\infty +t\alpha u_\infty)-\frac{\partial F}{\partial p}(q_\infty, p_\infty)\right)u_\infty dt  \leq \alpha.\eta= \frac{\alpha\varepsilon}{3}\hfill
 \end{matrix}$$
 and this is a contradiction.
 
 \enddemo
 We deduce from lemma \ref{Roc} that the sequence $(\frac{\partial H_k}{\partial p})$ $C^0$-converges   on any compact subset of $T^*M$.

Let us recall some well-known fact concerning Tonelli Hamiltonians that are proved for example in \cite{Fa1}.  To any Tonelli Hamiltonian $H$ (resp. $H_k$) we can associate a Lagrangian $L: TM\rightarrow \R$ (resp. $L_k: TM\rightarrow \R$) that is defined by:
 $$\forall (q, v)\in TM, L(q, v)=\sup_{p\in T_q^*M}\left( p.v-H(q,p)\right).$$
Then $L$ and $L_k$  are $C^2$, $C^2$-convex and superlinear in the fiber direction.\\
If we denote by $\Lc_q: T^*_qM\rightarrow T_qM$ the Legendre map $p\rightarrow \frac{\partial H}{\partial p}(q,p)$ (that is a $C^1$-diffeomorphism), we have the equivalent formula:
 $$L(q, v)=\Lc_q^{-1}(v).v-H(q, \Lc_q^{-1}(v)).$$
 Let us recall that the Euler-Lagrange flow  associated to $L$ is $f_t^L=\Lc\circ \varphi_t^H\circ \Lc^{-1}$ (see the notation $\Lc$ just below) and that if $(q_t, p_t)$ is an orbit for the Hamiltonian $H$, the corresponding orbit for $L$ is $(q_t, \dot q_t)$.

  We know that $\Lc^{H_k}: (q, p)\mapsto (q, \frac{\partial H_k}{\partial p}(q,p))=(q, \Lc_q^{H_k}(p))$ $C^0$-converges  on the compact subsets of $T^*M$ to  $\Lc : (q, p)\mapsto (q, \Lc_q(p))$. The following lemma implies that $((\Lc^{H_k})^{-1})$ $C^0$-converges    to $\Lc^{-1}$ on any compact subset of $TM$ and then that $(L_k)$ $C^0$-converges   to $L$ on any compact subset of $TM$.
 \begin{lemma}\label{homeo}
 Let $(h_k)$ be a sequence of homeomorphisms from $T^*M$ to $TM$ that $C^0$-converges    on any compact subset of $T^*M$ to a homeomorphism $h$. Then the sequence $(h_k^{-1})$ $C^0$-converges  on any compact subset of $TM$ to $h^{-1}$.
 \end{lemma}
  \noindent{\bf Proof of lemma \ref{homeo}.} If not, by possibly extracting   a subsequence,  we can find a   sequence  $(y_k)$ in $TM$ that converge to $y_\infty\in TM$ such that:
  $$\forall k, d(h_k^{-1}(y_k), h^{-1}(y_k))\geq 2\varepsilon>0.$$
  As $h$ is a homeomorphism, there exists $R_2>R_1>0$ such that:
  $$\forall x\in T^*M, d(h^{-1}(y_\infty), x)\geq R_2\Rightarrow d(y_\infty, h(x))\geq 4\varepsilon$$
  and 
  $$\forall x\in T^*M, d(h^{-1}(y_\infty), x)\leq R_1\Rightarrow d(y_\infty, h(x))<  \varepsilon.$$
  As $(h_k)$ $C^0$-converges   on any compact subset of $T^*M$, for $k$ large enough, we have: $h_k(B(h^{-1}(y_\infty), R_1))\subset B(y_\infty, \varepsilon)$ and\\
   $h_k(B(h^{-1}(y_\infty), 2R_2)\backslash B(h^{-1}(y_\infty), R_2))\subset TM\backslash B(y_\infty, 2\varepsilon)$. \\
  Hence $h_k^{-1}(B(y_\infty, \varepsilon))$ is a connected subset of $T^*M$ that meets $B(h^{-1}(y_\infty),R_1)$ but doesn't meet $B(h^{-1}(y_\infty), 2R_2)\backslash B(h^{-1}(y_\infty), R_2)$. This implies that    $h_k^{-1}(B(y_\infty, \varepsilon))\subset B(h^{-1}(y_\infty), R_2)$ .
  
    As the sequence $(y_k)$ converges to $y_\infty$, we may assume that $d(y_k, y_\infty)<\varepsilon$. Hence $d(h^{-1}(y_\infty), h_k^{-1}(y_k))\leq R_2$, the sequence $(h_k^{-1}(y_k))$ is bounded, we can extract a subsequence that converges to $x_\infty$. Then we have
    \begin{enumerate}
    \item[$\bullet$]  $d(x_\infty, h^{-1}(y_\infty))\geq 2\varepsilon$ by taking the limit in the first inequality;
    \item[$\bullet$] $\displaystyle{d(h(x_\infty), y_\infty)=\lim_{k\rightarrow +\infty} d(h_k(h_k^{-1}(y_k)), y_k)=0}$.
    \end{enumerate}
  We have found two   points $x_\infty\not=h^{-1}(y_\infty)$ that have the same image by the homeomorphism $h$. This is impossible.
  
  \enddemo
  Let us fix $(q, p)\in T^*M$. By Weierstrass theorem (see for example \cite{Fa1}), there exists $\tau>0$  such that the arc $\gamma: [0, \tau]\rightarrow M$ defined by $\gamma(t)=\pi\circ \varphi^H_t(q,p)$ is action strictly  minimizing for the Lagrangian $L$ , i.e. such that for every other absolutely continuous arc $\eta:[0, \tau]\rightarrow M$ that has the same endpoints as $\gamma$,  we have:
  $$A_L(\gamma)=\int_0^\tau L(\gamma(t), \dot\gamma(t))dt < A_L(\eta)=\int_0^\tau L(\eta(t), \dot\eta(t))dt.$$
  For any $k\in \N$, by Tonelli theorem (see \cite{Fa1}), there exists an arc $\gamma_k:[0, \tau ]\rightarrow M$ such that $\gamma_k(0)=q$ and $\gamma_k(\tau)=\gamma(\tau)$ that minimizes (non necessarily strictly)  the action for the Lagrangian $L_k$ among all the absolutely continuous arcs that have the same endpoints as $\gamma_k$. Then there exists a unique $p_k\in T_qM$ such that: $\gamma_k(t)=\pi\circ \varphi_t^{H_k}(q, p_k)$ and we have   $\dot\gamma_k(t)=\frac{\partial H_k}{\partial p}(\varphi_t^{H_k}(q, p_k))$. We define: $\displaystyle{m_\tau=\sup_{t\in[0, \tau]}|L(\gamma(t), \dot\gamma(t))}|$. 
  We will prove
  \begin{lemma}\label{Kompact}
  There exists a compact subset $K_2=\{ H\circ \Lc^{-1} \leq C+2\}$ of $TM$ that contains all the $(\gamma_k(t), \dot\gamma_k(t))$ and $(\gamma(t), \dot\gamma(t))$ for $t\in [0, \tau]$ and we have
  $$\lim_{k\rightarrow +\infty} A_L(\gamma_k)=A_L(\gamma).$$
  \end{lemma}
  \noindent{\bf Proof of lemma \ref{Kompact}.}\\
  As $(L_k)$ $C^0$-converges   to $L$ on the compact subsets of $TM$, for $k$ large enough, we have: $A_{L_k}(\gamma)\leq \tau (m_\tau+1)$ and then $A_{L_k}(\gamma_k)\leq \tau(m_\tau+1)$ because $\gamma_k$ is minimizing for $L_k$.\\
Because $L$ is superlinear and convex in the fiber direction, there exists $R>0$ such that 
$$\forall (q,v)\in TM, \| v\|\geq R \Rightarrow L(q, v)> M_\tau= \sup\left\{  m_\tau+3, \sup_{q\in M}|L(q, 0)|+1\right\}.$$
Because $L_k$ is convex in the fiber direction, we have for $\| v\|\geq R$
$$L_k(q, v)\geq L_k(q, 0)+\left(L_k(q,  R\frac{v}{\| v\|})-L_k(q, 0)\right) \frac{\| v\| -R}{R}.$$
Because $(L_k)$ $C^0$-converges to $L$ on the compact subsets of $TM$, for $k$ large enough, we have
$$\forall (q, v), L_k(q,0)\geq -M_\tau\quad{\rm and}\quad L_k(q, R\frac{v}{\| v\|})-L_k(q, 0)>1.$$
We deduce for $k$ large enough and $\| v\|\geq R$
$$L_k(q, v)\geq   \frac{\| v\| -R}{R}-M_\tau$$
and then for a certain $R^*>R$ , we have
$$\forall v, \| v\|\geq R^*\Rightarrow L_k(q, v)\geq M_\tau.$$

   For $k$ large enough, we have $A_{L_k}(\gamma_k)\leq \tau (m_\tau +1)$ and then for a $t_k\in[0, \tau]$ we have $\|\dot\gamma_k(t_k)\|\leq R^*$.  We introduce the constant $\displaystyle{C=\sup_{\| v\|\leq R^*}H(  \Lc^{-1}(q, v))}$ and the compact subsets $K_1=\{ H\circ \Lc^{-1} \leq C+1\}$,  $K_2=\{ H\circ \Lc^{-1} \leq C+2\}$ and $K_3=\{ H\circ \Lc^{-1} \leq C+3\}$  of $TM$. As $H$ is convex on the fiber directions the intersection of $K_i$ with a fiber is a topological ball. As $(H_k\circ ( \Lc^{H_k})^{-1})$ $C^0$-converges to $ H\circ \Lc^{-1}$ on $K_2$, for $K$ large enough we have
   $$\forall (q, v)\in K_3\backslash K_2, H_k\circ  ( \Lc^{H_k})^{-1}(q,v)> C+\frac{3}{2}  $$
  and $$\forall q, \exists (q,v)\in K_1, H_k\circ ( \Lc^{H_k})^{-1}(q,v)<C+\frac{3}{2}.$$
  As $H_k$ is convex in the fiber direction, the set $\{ H_k\circ  ( \Lc^{H_k})^{-1}\leq C+1\}$ is connected and doesn't meet the boundary of $K_2$;  hence $\Ec_k=\{ H_k\circ  ( \Lc^{H_k})^{-1}\leq C+1\}\subset K_2$.

 The set $\Ec_k$ is invariant by the Euler-Lagrange flow of $L_k$ (the energy is invariant); we know that  $\|\dot\gamma_k(t_k)\|\leq R^*$ and $H\circ \Lc^{-1}(\gamma_k(t_k), \dot\gamma_k(t_k))\leq C$, hence for $k$ large enough: $H_k\circ (\Lc^{H_k})^{-1}(\gamma_k(t_k), \dot\gamma_k(t_k))\leq C+1$ and then $(\gamma_k(t_k), \dot\gamma_k(t_k))\in \Ec_k\subset K_2$ and
 $$\forall t\in [0, \tau],( \gamma_k(t), \dot\gamma_k(t))\in\Ec_k\subset K_2.$$
  
 We introduce the notation $\displaystyle{\varepsilon_k=\sup_{(q,v)\in K_2}|L_k(q, v)-L(q,v)|}$. We have then: 

   $A_{L }(\gamma_k)\leq A_{L_k}(\gamma_k)+\tau\varepsilon_k \leq A_{L_k}(\gamma)+\tau\varepsilon_k\leq A_L(\gamma)+2\tau\varepsilon_k$. We deduce that $\displaystyle{\lim_{k\rightarrow +\infty}A_L(\gamma_k)=A_L(\gamma)}$. Let us notice that we obtain the same result by replacing $\tau$ by any smaller $\tau$.\enddemo
   
   \begin{remk}
   Let us explain the link between  the  previous proof and some proofs that are given by P.~Bernard in \cite{Ber1}. He introduces a notion of ``uniform families of Tonelli Hamiltonians'' that satisfy in particular
   \begin{enumerate}
   \item[(i)] there exist $h_0, h_1:\R^+\rightarrow \R^+$ that are superlinear such that for every Hamiltonian $H$ of the family
   $$h_0(\| p\|)\leq H(q, p)\leq h_1(\| p\|);$$
   \item[(ii)] there exists an increasing function $k:\R^+\rightarrow \R_+$ such that if $|t|\leq 1$ and $(\varphi_s)$ is the flow associated to a Hamiltonian of the family, then
   $$\varphi_t(\{ \| p\|\leq r\} \subset \{ \| p\|\leq k(r)\}\leq T^*M.$$
   \end{enumerate}
   It is easy to prove that if the sequence $(H_n)$ of Tonelli Hamiltonians $C^0$-converges on the compact subsets to a Tonelli Hamiltonian $H$, then the family $\{ H_n\}$ satisfy (i) and (ii). With these conditions, P.~Bernard proves in \cite{Ber1} (see (B5))a result that implies the first part of lemma \ref{Kompact}: the minimizers in time $[0, t]$ of all the Hamiltonians of the family have a derivative that is uniformly bounded.
   \end{remk}

   Let us now give the proof of proposition \ref{Tonelli} that is
 \begin{propo}{\bf 1}   Let $(H_k)_{k\in\N}$, $H$ be $C^2$ Tonelli Hamiltonians defined on $T^*M$ such that $(H_k)$ $C^0$-converges   on the compact subsets of $T^*M$ to $H$. 
 
 Then, for every $(q,p)\in T^*M$, there exists $T>0$ and $(q_k, p_k)_{k\in \N}\in T^*M$  such that a subsequence of the  sequence of arcs $(Q_k, P_k)=(\varphi_t^{H_k}(q_k, p_k))_{t\in [0, T]}$ converges to  $(Q, P)=(\varphi_t^H(q, p))_{t\in[0, T]}$ in the following sense (here we are in a chart because all is local and we denote the subsequence in the same way we denoted the initial sequence):
 \begin{enumerate}
 \item[$\bullet$] $(Q_k)$ $C^0$-converges to $Q$;
 \item[$\bullet$] $(P_k)$ converges to $P$ in the $L^2$ sense;
 \item[$\bullet$] the $(Q_k, P_k)$ are uniformly $C^0$ bounded.
 \end{enumerate} 
 \end{propo}

  Lemma \ref{Kompact} implies that the family $(\gamma_k)$ is equicontinuous (see \cite{Fa1}) and then $C^0$-relatively compact by Ascoli theorem. A classical result (see corollary 3.2.3. in \cite{Fa1}) asserts that $A_L$ is lower semi-continuous for the $C^0$-topology. This implies that every limit point of the sequence $(\gamma_k)$ is $\gamma$ and then that $(\gamma_k)$ $C^0$-converges to $\gamma$. Assuming that $k$ is large enough, we can then assume that all the $\gamma_k$'s and $\gamma$ are in a same chart $B(\gamma(0), R)$. We have proved too that the $\dot\gamma_k$ are uniformly bounded.\\
  
 Let us now prove that   $(\dot\gamma_k)$   converges   to $\dot\gamma$ for $\| .\|_2$. The method that we use was suggested by Patrick Bernard. \\
 We use the notation $p(t)=(\Lc^{H})^{-1}(\gamma(t), \dot\gamma(t))$. Using the method of characteristics (see for example \cite{Fa1}), it is easy to build in some neighbourhood of $\gamma(0)$ a $C^2$ local solution $u_t: B(\gamma(0), r)\rightarrow \R$ of the Hamilton-Jacobi equation such that $du(\gamma(0))=p(0)$ (with $t\in [-\tau, \tau]$ with an possibly smaller $\tau$); we have then
 $$\forall q\in B(\gamma(0), r), \frac{\partial u_t}{\partial t} (q)+H(q, du_t(q))=0.$$
 As $\gamma$ is a solution of the Euler-Lagrange equations for $L$, note that we have the equality in the Young inequality
 $$H(\gamma(t), du_t(\gamma (t)))+L(\gamma(t), \dot\gamma(t))=du_t(\gamma(t))\dot\gamma(t);$$
 integrating the two previous equalities, we deduce that for every $t\in(0, \tau]$, we have 
 $$A_L(\gamma_{|[0, t]})=u_t(\gamma(t))-u_0(\gamma(0)).$$
 Now, Young inequality gives
 $$L(q, v)+H(q, du_t(q))-du_t(q)v\geq 0.$$
The second derivative with respect to $v$ of the above function   is bounded from below on every compact subset. Hence there exists $C>0$ such that
$$\forall (q, v)\in K_2 \cap \pi^{-1}(B(\gamma(0), r)), L(q, v)+H(q, du_t(q))-du_t(q)v\geq C\| v-(\Lc^H)^{-1}(du_t(q))\|^2;$$
i.e. 
 $$\forall (q, v)\in K_2 \cap \pi^{-1}(B(\gamma(0), r)), L(q, v)-\frac{\partial u_t}{\partial t}(q)-du_t(q)v\geq C\| v-(\Lc^H)^{-1}(du_t(q))\|^2.$$
 We take $(q, v)=(\gamma_k(t), \dot\gamma_k(t))$ and integrate between $0$ and 
$ \tau$
 $$ A_L(\gamma_k)-(u_\tau(\gamma_k(\tau))-u_0(\gamma_k(0))\geq C\int_0^\tau \|\dot\gamma_k (t)-(\Lc^H)^{-1}(du_t(\gamma_k(t)))\|^2dt.$$
 As $(\gamma_k)$ uniformly converge to $\gamma$, the left hand term converges to zero when $k$ tends to $+\infty$ and $t\mapsto (\Lc^H)^{-1}(du_t(\gamma_k(t)))$ uniformly converges to $\dot\gamma(t)=(\Lc^H)^{-1}(du_t(\gamma(t)))$. We deduce that 
 $$\lim_{k\rightarrow +\infty}\| \dot\gamma_k-\dot\gamma\|_2=0.$$

 \begin{remk} As suggested by the referee, we suggest another possible proof of the convergence for $\| .\|_2$. We don't assume that $H$ has any polynomial growth in the fiber direction.  But by using \cite{ConItu1}, we can easily modify $H$ outside a large enough compact to obtain a $H$ with quadratic growth. In this case, it is possible to use proposition 3.10 p 123 of the book \cite{ButGiaHil} by G.~Buttazzo, M.~Giaquinta \& S.~Hildebrandt to obtain the convergence for $\|.\|_2$.
 
 \end{remk}

  The Lagrangians $L_k$, $L$ are $C^2$-convex in the fiber direction  and the sequence $(L_k)$ $C^0$-converges    on any compact subset of $TM$.  Hence, by lemma \ref{Roc},  the sequence $(\frac{\partial L_k}{\partial v})$ $C^0$-converges  on any compact subset of $TM$ to $\frac{\partial L}{\partial v}$. We deduce that $(\varphi_t^{H_k}(q, p_k))_{t\in[0, T]}=(\gamma_k(t), \frac{\partial L_k}{\partial v}(\gamma_k(t), \dot\gamma_k(t)))_{t\in [0, T]}$ is uniformly bounded  and converges  to $(\varphi_t^H(q, p))_{t\in[0, T]}=(\gamma(t), \frac{\partial L}{\partial v}(\gamma(t), \dot\gamma(t)))_{t\in[0, T]}$ for the distance that is the sum (we are in a chart) of the uniform distance of $M$ and the $\| .\|_2$ distance in the fiber direction. Indeed, if we use the notation $\varepsilon_k=\sup_{K_2}\|\frac{\partial L}{\partial v}-\frac{\partial L_k}{\partial v}\|$ and if $\ell$ is a Lipschitz constant for $\frac{\partial L}{\partial v}$ on $K_2$, we have (in a chart)
  $$\begin{matrix}   \| \frac{\partial L_k}{\partial v}(\gamma_k , \dot\gamma_k)-&\frac{\partial L}{\partial v}(\gamma , \dot\gamma ) \|_2\hfill\\
  &\leq \| \frac{\partial L_k}{\partial v}(\gamma_k , \dot\gamma_k)-\frac{\partial L}{\partial v}(\gamma_k, \dot\gamma_k) \|_2+\| \frac{\partial L} {\partial v}(\gamma_k , \dot\gamma_k) -\frac{\partial L}{\partial v}(\gamma, \dot\gamma) \|_2\\
  &\leq \sqrt{\tau}\varepsilon_k+\ell\left(  \sqrt{\tau}\|\gamma_k-\gamma\|_{C^0}+\|\dot\gamma_k-\dot\gamma\|_2  \right)  \hfill\end{matrix}$$

 \enddemo
 \section{Proof of propositions \ref{clef} and \ref{2clef}}\label{Sec4}
 \subsection{Proof of proposition \ref{clef}}
 We recall proposition \ref{clef}
 \begin{propo}{\bf 2}
  Let $(F_k)$ be a sequence of $C^2$ Hamiltonians defined on $T^*M$ and let $F$ be a $C^2$ Hamiltonian such that:\\
   for every $(q,p)\in T^*M$, there exists $T>0$ and $(q_k, p_k)_{k\in \N}\in T^*M$  such that the sequence of arcs $(Q_k, P_k)=(\varphi_t^{F_k}(q_k, p_k))_{t\in [0, T]}$ converges to  $(Q, P)=(\varphi_t^F(q, p))_{t\in[0, T]}$ in the following sense:
 \begin{enumerate}
 \item[$\bullet$] $(Q_k)$ $C^0$-converges to $Q$;
 \item[$\bullet$] $(P_k)$ converges to $P$ in the $L^2$ sense;
 \item[$\bullet$] the arcs $(Q_k, P_k)$ are uniformly $C^0$ bounded.
 \end{enumerate} 
Let $(G_k)$ be a sequence of $C^2$ Hamiltonians of $T^*M$ that $C^0$-converges on the compact subsets of $T^*M$ to a $C^2$ Hamiltonian $G$ and that is such that the sequence $(\{ F_k, G_k\})$ $C^0$-converges on the compact subsets of $T^*M$ to some $\Hc\in C^0(T^*M, \R)$.\\
Then $\{ F, G\} =\Hc$. 
\end{propo}
 \noindent Let $(F_k)$, $F$, $(G_k)$, $G$ as in the hypotheses of proposition \ref{clef}. \\
 Let $(q, p)\in T^*M$ be a point. We want to prove that $\{F, G\}(q, p)=\Hc(q,p)$. \\
 We associate $T$, $(q_k, p_k)$, $(Q, P)$, $(Q_k, P_k)$ to $(q, p)$ as above.
 As $(P_k)$ converges to $P$ in the $L^2$ sense,  $(P_k)$ has a subsequence (also denoted by $(P_k)$) that converges almost everywhere to $P$. Because $(G_k)$ converge to $G$ on the compact subsets  and the $P_k$ are uniformly bounded, we deduce that  for almost every $s<t$ in $[0, T]$, we have
 $$\begin{matrix} G(Q(t), P(t))-G(Q(s), P(s))&=\displaystyle{\lim_{k\rightarrow +\infty} (G_k(Q_k(t), P_k(t))-G_k(Q_k(s), P_k(s))}\\
 &=\displaystyle{\lim_{k\rightarrow +\infty}\int_{s}^{t} \{ G_k, F_k\} \circ (Q_k(\sigma), P_k(\sigma))d\sigma.}\hfil\end{matrix}$$
 Note that 
 $$\begin{matrix}&\left|\int_{s}^{t}\  (\{ G_k, F_k\} \circ (Q_k(\sigma), P_k(\sigma))+\Hc(Q(\sigma), P(\sigma)))d\sigma\right|\hfill\\
 &\leq \left| \int_s^t(\{ G_k, F_k\}+\Hc)(Q_k(\sigma), P_k(\sigma))d\sigma\right|+\left| \int_s^t(\Hc(Q_k(\sigma), P_k(\sigma))-\Hc(Q(\sigma), P(\sigma)))d\sigma\right|.\end{matrix}
 $$
 The first term tends to $0$ because of the convergence of $(\{ F_k, G_k\})$ to $\Hc$ on the compact subsets and the second one tends to $0$ by the Lebesgue dominated convergence theorem.
 
 We finally obtain for almost every $s$, $t$
 
 $$G(Q(t), P(t))-G(Q(s), P(s))=-\int_s^t \Hc(Q(\sigma), P(\sigma))d\sigma$$
 and by continuity the result for all $s$, $t$.
Differentiating at $0$ with respect to $t$, we obtain
 $$\{ G, F\} (q, p)=-\Hc(q,p).$$
 \enddemo
 \subsection{Proof of proposition \ref{2clef}}
 The idea is exactly the same as for proposition \ref{clef}. the proof is simpler because we assume a $C^0$ convergence. Let us recall proposition \ref{2clef}.
 \begin{propo}{\bf 3}
 Let $(N, \omega)$ be a (non necessarily compact) symplectic manifold. Let $(F_k)$ be a sequence of $C^2$ Hamiltonians defined on $N$ and let $F$ be a $C^2$ Hamiltonian such that:\\
   for all $x\in N$, there exists  $T>0$ and a sequence $(x_k)\in N$ such that the sequence of arcs of orbit $(\varphi_t^{F_k}(x_k))_{t\in [0, T]}$ $C^0$-converges to $(\varphi_t^{F}(x))_{t\in [0, T]}$.\\
Let $(G_k)$ be a sequence of $C^2$ Hamiltonians of $N$ that $C^0$-converges on the compact subsets of $N$ to a $C^2$ Hamiltonian $G$ and that is such that the sequence $(\{ F_k, G_k\})$ $C^0$-converges to some $\Hc\in C^0(N, \R)$.\\
Then $\{ F, G\} =\Hc$. 

 \end{propo}
 Because of the $C^0$-convergence of  $(\varphi_t^{F_k}(x_k))_{t\in [0, T]}$ to $(\varphi_t^{F}(x))_{t\in [0, T]}$, we can choose a compact subset $K$ of $N$ that contains all the $\varphi_t^{F_k}(x_k)$. As $(G_k)$ $C^0$-converges to $G$ on $K$, we have
 $$G(\varphi_T^F(x))-G(x)=\lim_{k\rightarrow +\infty} G_k(\varphi_T^{F_k}(x))-G(x)=\lim_{k\rightarrow +\infty}�\int_0^T\{G_k, F_k\}(\varphi_t^{F_k}(x))dt.$$
 As $\{ G_k, F_k\}$ $C^0$-converges to $-\Hc$ on $K$ and $(\varphi_t^{F_k}(x_k))_{t\in [0, T]}$ $C^0$-converges to $(\varphi_t^{F}(x))_{t\in [0, T]}$, we have
 $$\lim_{k\rightarrow +\infty} �\int_0^T\{G_k, F_k\}(\varphi_t^{F_k}(x))dt=-\int_0^T\Hc(\varphi_t^F(x))dt.$$
 Differentiating $G(\varphi_T^F(x))-G(x)=-\int_0^T\Hc(\varphi_t^F(x))dt$ with respect to $T$, we obtain $\{ G, F\} (x)=-\Hc(x)$.\enddemo


\begin{thebibliography}{cc}

\bibitem {Ber1} P.~Bernard. \emph{The dynamics of pseudographs in convex Hamiltonian systems.}  J. Amer. Math. Soc.  21  (2008),  no. 3, 615--669.  

\bibitem{Buh1} L.~Buhovsky, The 2/3-convergence rate for the Poisson bracket. Geom. Funct. Anal. {\bf 19-6}  (2010),  1620-1649

\bibitem{ButGiaHil} G.~Buttazzo, M.~Giaquinta \& S.~Hildebrandt, {\em One-dimensional variational problems. An introduction} Oxford Lecture Ser. Math. Appl. {\bf 15}, The Clarendon Press, Oxford university Press, New York, 1998.

\bibitem{CarVit}  F.~Cardin \& C.~Viterbo, {\em Commuting Hamiltonians and Hamilton-Jacobi multi-time equations},  Duke Math. J. {\bf 144-2}  (2008),   235-284. 

\bibitem {ConItu1}  G.~Contreras \&  R.~Iturriaga,  {\em Convex
Hamiltonians without conjugate points}  Ergodic Theory Dynam. Systems  {\bf 19 }
,  no. 4,  (1999)901--952. 

\bibitem{CIPP1} G.~Contreras, R.~Iturriaga, G.~P.~ Paternain \& M.~Paternain, {\em The Palais-Smale condition and Ma\~n\'e's critical values},  Ann. Henri Poincar\'e {\bf 1-1}  (2000),  655-684.

\bibitem{Fa1} A.~Fathi, {\em Weak KAM theorems in Lagrangian dynamics},
book in preparation. 

\bibitem{HubWes} J.~H.~Hubbard \& B.~H.~ West, {\sl Equations diff\'erentielles et syst\`emes dynamiques},  Cassini, Paris, 1999 xiv+416pp
\bibitem{Hum1}  V.~Humili\`ere, {\em Continuit\'e en topologie symplectique}, PhD thesis, july 2008

\bibitem{Hum2}  V.~Humili\`ere, {\em Pseudo-representations},  Comment. Math. Helv. {\bf 84-3} (2009),  571-585.


\bibitem{Roc1}  R.~T.~Rockafellar,   Convex analysis. Princeton Mathematical Series, No. 28 Princeton University Press, Princeton, N.J. 1970 xviii+451 pp


\end{thebibliography}
\end{document}